\documentclass[12pt]{article}

\title{A la Fock-Goncharov coordinates for PU(2,1)}

\usepackage{amsmath,amssymb,amsthm,graphicx,pstricks,xypic,comment}

\begin{document}
\author{Julien March\'e, Pierre Will\\
        Institut de Math{\'e}matiques\\
        Universit{\'e} Pierre et Marie Curie\\
        4, place Jussieu\\
        F-75252 Paris Cedex 05\\ 
       email :  \texttt{marche,will@math.jussieu.fr}}

\newcommand{\la}{\langle}
\newcommand{\ra}{\rangle}
\newcommand{\HdC}{\bold H^{2}_{\mathbb C}}
\newcommand{\HdR}{\bold H^{2}_{\mathbb R}}
\newcommand{\HuC}{\bold H^{1}_{\mathbb C}}
\newcommand{\HnC}{\bold H^{n}_{\mathbb C}}
\newcommand{\Cdu}{\mathbb C^{2,1}}
\newcommand{\Ct}{\mathbb C^{3}}
\newcommand{\C}{\mathbb C}
\newcommand{\Mn}{\mbox{M}_n\left(\mathbb C\right)}
\newcommand{\R}{\mathbb R}
\newcommand{\A}{\mathbb A}
\newcommand{\cC}{\mathfrak C}
\newcommand{\cCp}{\mathfrak C_{\mbox{p}}}
\newcommand{\cT}{\mathfrak T}
\newcommand{\rR}{\mathfrak R}
\newcommand{\hH}{\mathfrak H}
\newcommand{\Td}{T^{2}_{\left(1,1\right)}}
\newcommand{\T}{{\bf T}}
\newcommand{\PSL}{\mbox{\rm PSL(2,$\R$)}}
\newcommand{\Pu}{\mbox{\rm PU(2,1)}}
\newcommand{\X}{\mbox{\textbf{X}}}
\newcommand{\tr}{\mbox{\rm tr}}
\newcommand{\M}{{\rm M}}
\newcommand{\N}{{\rm N}}
\newcommand{\B}{\mathcal{B}}
\newcommand{\n}{\noindent}
\newcommand{\p}{{\bf p}}
\newcommand{\q}{{\bf q}}
\renewcommand{\Re}{\mbox{\rm Re}\,}
\renewcommand{\Im}{\mbox{\rm Im}\,}
\renewcommand{\leq}{\leqslant}
\renewcommand{\geq}{\geqslant}
\newcommand{\mat}[1]{\begin{bmatrix} #1 \end{bmatrix}}
\newcommand{\ba}[1]{\overline{#1}}
\newcommand{\bc}{{\bf c}}
\newcommand{\bd}{{\bf d}}
\newcommand{\bp}{{\bf p}}
\newcommand{\bm}{{\bf m}}
\newcommand{\bn}{{\bf n}}
\newcommand{\boX}{\mathcal{X}}
\newcommand{\boM}{\mathcal{M}}
\renewcommand{\phi}{\varphi}

\maketitle
\newtheorem{theo}{Theorem}
\newtheorem{lem}{Lemma}
\newtheorem*{lem*}{Lemma}
\newtheorem{coro}{Corollary}
\newtheorem{prop}{Proposition}
\theoremstyle{remark}
\newtheorem{rem}{Remark}
\newtheorem{ex}{Example}
\theoremstyle{definition}
\newtheorem{defi}{Definition}
\abstract{In this work, we describe a set of coordinates on the PU(2,1)-representation variety of the fundamental group of an oriented punctured surface $\Sigma$ with negative Euler characteristic. The main technical tool we use is a set of geometric invariants of a triple of flags in the complex hyperpolic plane $\HdC$. We establish a bijection between a set of decorations of an ideal triangulation of $\Sigma$ and a subset of the PU(2,1)-representation variety of $\pi_1(\Sigma)$.}

\section{Introduction}

In their work \cite{FockGon3}, Fock and Goncharov have described a coordinate system on the 
representation variety of the fundamental group of a punctured surface $\Sigma$ in a split semi-simple 
real Lie group. The typical example is PSL($n$,$\R$). When $n$ equals 2 (resp. 3), they identify the Teichm\"uŸller 
space of $\Sigma$ (resp. the moduli space of convex projective structures on $\Sigma$) within the representation 
variety. The main goal of our work is to describe an analogous coordinate system for representations in PU(2,1), 
which is the group of holomorphic isometries of the complex hyperbolic plane $\HdC$. Note that PU(2,1) is not split 
and thus does not belong to the family studied by Fock and Goncharov. The preprint \cite{FockGon2} in which the cases of 
PSL($2$,$\R$) and PSL($3$,$\R$) are dealt with separately has been our main source of inspiration (see
also \cite{FockGon4}).

Throughout this article, we will use the following notation. Let $\Sigma_{g,p}$ be a genus 
$g$ surface with $p$ punctures $x_1,\ldots,x_p$, assuming $p>0$ and $2-2g-p<0$. We denote by $\pi_{g,p}$ its fundamental group and 
use the following standard presentation where the $c_i$'s are homotopy classes of curves 
enclosing the $x_i$'s:

$$\pi_{g,p}=\la a_1,b_1,\ldots 
a_g,b_g,c_1,\ldots,c_p\vert\prod_{i=1}^g\lbrack
a_g,b_g\rbrack\prod_{j=1}^p c_j\ra.$$

We will call \textit{flag} a pair $(C,p)$ where $C$ is a complex line of $\HdC$ (see definition \ref{complexline}) and 
$p$ is a boundary point of $C$. Our goal is to parametrize the variety

$$\mathfrak{R}_{g,p}=\left\{\rho,\mathcal{F}\right\}/PU(2,1)$$
where 
\begin{itemize}
\item
$\rho$ is a representation of $\pi_{g,p}$ in PU(2,1)
\item
$F=(F_1,\ldots,F_p)$ is a $p$-tuple of flags such that $\rho(c_i)$ stabilizes $F_i$. 
\item
The group PU(2,1) acts on $\rho$ by conjugation and on $F$ by isometries.
\end{itemize}

\begin{rem}
In \cite{FockGon2,FockGon3}, the authors use an alternative
description of $\mathfrak{R}_{g,p}$ which is equivalent to the above
one but appears to  
be more efficient for certain aspects. We recall it briefly for later
use. Let $\widehat{\Sigma}_{g,p}$ be the universal covering of
$\Sigma_{g,p}$. The surface $\widehat{\Sigma}_{g,p}$ may be seen as a
topological disk with an action of $\pi_{g,p}$ and an invariant family
$X$ of boundary points projecting onto the $x_i$'s.  
The space $\mathfrak{R}_{g,p}$ is in bijection with the space of couples
$(\rho,F)$ up to conjugation, where $\rho$ is a representation of
$\pi_{g,p}$ in PU(2,1) and $F$ is an equivariant map from $X$ to the
space of flags, that is, for all $g$ in $\pi_{g,p}$ and $x$ in $X$ one
has $F(g.x)=\rho(g).F(x)$. 

Let us give some rough indications about the equivalence of the two
definitions. The curve $c_i$ determines a prefered lift of $x_i$ in
$X$ which we denote by $\widehat{x_i}$. Given  
an equivariant map $F$, we just set
$F_i=F(\widehat{x_i})$. Reciprocally, if we have a $p$-tuple
$(F_1,\ldots,F_p)$ of flags, we construct a map $F$ by setting
$F(\widehat{x_i})=F_i$ and extend it by the equivariance property. 

\end{rem}

To build the coordinate system, we start from an ideal triangulation
$T$ of $\Sigma_{g,p}$ (see definition \ref{triangul}). Each triangle
$\Delta$ of $T$ lifts to $\widehat{\Sigma}_{g,p}$ as a triangle whose
vertices are denoted by $x,y,z\in X$. Given a pair $(\rho,F)$, the
triple of flags $\left(F(x),F(y),F(z)\right)$ is well defined up to
isometry. In definition \ref{genflag}, we introduce a notion of
genericity for a triple of flags. We will say that a couple $(\rho,F)$
is {\it generic with respect to $T$} if the triple of flags associated
to any triangle of $T$ is generic. 
We will denote by $\mathfrak{R}_{g,p}^{T}$ the subset of
$\mathfrak{R}_{g,p}$ containing those class of pairs $(\rho,F)$ that
are generic with respect to $T$.

Next, we associate to $\Delta$ a family of invariants which parametrizes generic triples of flags up to isometry.
The definition and study of these invariants is a crucial point of the
article. To represent a geometric configuration,  
the invariants associated to adjacent triangles must satisfy compatibility
relations. We will call {\it decoration} of a  
triangulation $T$ the following data: a family of invariants for each
triangle of $T$ such that  the compatibility conditions are satisfied
(see definition \ref{defdeco}). 

The main result of the article is the following

\begin{theo}\label{bij}
Let $T$ be an ideal triangulation of  $\Sigma_{g,p}$. There is a bijection between 
$\mathfrak{R}_{g,p}^{T}$ and $\boX(T)$, the set of decorations of $T$.
\end{theo}

Note that some isometries of $\HdC$ do not preserve any flag. These isometries are 
unipotent parabolic and are conjugate in PU(2,1) to non-vertical Heisenberg translations (see chap. 4 in \cite{Go}). 
As a consequence, $\mathfrak R_{g,p}$ do not contain all representations of $\pi_{g,p}$ in PU(2,1).

To any ideal triangulation, Fock and Goncharov associate a coordinate system on the representation 
variety. The transition from a triangulation to another may be done by a succession of elementary moves, the 
so-called flips, which allow them to forget about the initial choice of a triangulation. 
In the cases treated by Fock and Goncharov, the introduction of coordinate systems gives rise to a special class 
of representations called positive. By computing the coordinate
changes associated to the flips, they show first that the  
positivity of a representation is independent of the choice of
triangulation and second that the positive representations are  
discrete and faithful. These coordinates appear {\it a posteriori} to
be a quick and elegant way to study Teichm{\"u}ller spaces  
and their generalizations. Such a treatment of discreteness in the case of PU(2,1) 
seems to be still out of reach.

The study of representations of surface groups in PU(2,1) began in the eighties with Goldman and 
Toledo among others (see \cite{Go7,Tol}). However, many natural questions still do not have received a 
complete answer. Apart from a few general results about rigidity and flexibility (see \cite{Go7,PP1,Tol}) 
most of the results are dealing with examples or families of examples (see \cite{FK,GKL,Wi2}). 

Up to this day, no example of a PU(2,1)-moduli space of discrete and faithful representations of a
given surface group has been described. The only infinite group of finite type for which all the 
discrete and faithful representations in PU(2,1) are known is the modular group PSL(2,$\mathbb{Z}$) (see 
\cite{FJP}). In the case of closed surfaces, Parker and
Platis have described in \cite{PP} coordinates analogous to Fenchel-Nielsen
coordinates in the setting of PU(2,1).

In this article, we have chosen to introduce all the notions of complex hyperbolic geometry we are using.
 Some of the invariants we are dealing with are very classical. As an example, the invariant $\phi$ of a pair of complex 
 lines (see \ref{2lines}) is treated in \cite{Go} and the classification of triples of complex lines (see \ref{3lines}) 
 appears in \cite{Pra}. We decided to include the definitions and proofs about these invariants for the convenience of 
 the reader. Nevertheless, the invariants $m$ and $\delta$ (see definitions \ref{defm} and \ref{defdelta}) are specially adapted to 
 pairs of flags and do not appear elsewhere to our knowledge.

The article is organized as follows:

\begin{itemize}
\item
The section \ref{HdC} is devoted to the exposition of notions of complex 
hyperbolic geometry. We describe totally geodesic subspaces of the complex 
hyperbolic plane and introduce invariants of pairs and triples of complex lines.
\item
In section \ref{flaginv}, we describe the main technical tools which are the 
invariants $m$ and $\delta$. The main result of this section is the theorem \ref{classdrap} 
which classifies triples of flags up to isometry.
\item
In section \ref{stand}, we define the {\it standard configuration} of a flag and a complex line. Using 
the invariants described in the previous section, we provide two explicit matrices that are the elementary pieces necessary to construct the representations from the invariants. These matrices may be useful for numerical applications.
\item
The section \ref{deco} is devoted to the definition of the decoration space and to the proof 
of theorem \ref{bij}.
\item We prove in section \ref{resol} that the compatibility equations involved in the decoration space always
have solutions. The main tool is the lemma \ref{resou3}: it shows that once the $\Phi$
and $m$ invariants of a triple of flag are known, there exist generically 2 possible triples
of $\delta$ invariants, which correspond to the fixed points of an antiholomorphic isometry in the 
boundary of a disk. As a consequence of this lemma, we obtain in proposition \ref{cover} 
that $\mathfrak{R}_{g,p}^T$ is a $2^N$ cover of a simpler space denoted by $\mathcal{M}_{g,p}^T$, 
which is an auxilliary decoration space of the triangulation $T$, given by the $\Phi$ and $m$ invariants.

\item We give in section \ref{cusps} some indications about how to control the isometry type of the
images of the boundary curves by the representation constructed from a decorated triangulation of
$\Sigma_{g,p}$. We first deal with the case of an arbitrary punctured surface, and move then to the case
of the 1-punctured torus.
\end{itemize}

\section{Complex hyperbolic geometry\label{HdC}}
\subsection{Generalities}
Consider the hermitian form of signature $(2,1)$ in $\C^3$ given by the formula $\la v,w\ra
=v^T J\ba{w}$ where $J$ is the matrix given by 
$$J=\begin{bmatrix} 0&0&1\\ 0&1&0\\ 1&0&0\end{bmatrix}.$$ 

\n We define the following subsets of $\C^3$: 
\begin{eqnarray*}
V_0=\{ v\in \C^3\setminus\{ 0\}, \langle v,v\rangle =0\}\\
V_-=\{ v\in \C^3\setminus\{ 0\}, \langle v,v\rangle <0\}\\
V_+=\{ v\in \C^3\setminus\{ 0\}, \langle v,v\rangle >0\}
\end{eqnarray*}

\n Let ${\bf P}:\C^3\setminus\{0\}\to \C P^2$ be the canonical projection onto
the complex projective space. 

\begin{defi}
The complex hyperbolic plane $\HdC$ is the set ${\bf P}(V_-)$
equipped with the Bergman metric. 
\end{defi}
\n The boundary of $\HdC$ is ${\bf P}\left(V_0\right)$.
\n The distance function associated to the Bergman metric is given in
terms of Hermitian product by 
\begin{equation}\cosh^2\left(\dfrac{d(m,n)}{2}\right)=\dfrac{\la {\bf m},{\bf n}\ra\la{\bf n},
{\bf m}\ra}{\la{\bf m},{\bf m}\ra\la{\bf n},{\bf n}\ra},\label{distance}\end{equation}
where ${\bf m}$ and ${\bf n}$ are lifts of $m$ and $n$ to $\C^3$.
\n It follows from (\ref{distance}) that U(2,1), the unitary group
associated to $J$, acts on $\HdC$ by holomorphic isometries. The full
isometry group of $\HdC$ is generated by PU(2,1) and the complex
conjugation. The usual trichotomy of isometries for PSL(2,$\R$) holds here also: an isometry 
is elliptic if it has a fixed point inside $\HdC$, parabolic if it has a unique fixed point on 
$\partial\HdC$, and loxodromic if it has exactly two fixed points on $\partial\HdC$, and this 
exhausts all possibilities.

\subsection{Subspaces of $\HdC$.}
There are two types of maximal totally geodesic subspaces of $\HdC$, which are both of (real) 
dimension 2: complex lines and $\R$-planes. We give now a few indications about these. More details
may be found in \cite{Go}.

\subsubsection{Complex lines}
\begin{defi}\label{complexline}
We call \textit{complex line} in $\HdC$ the intersection with $\HdC$ of the projectivization of a 
2-dimensional subspace of $\C^3$ which intersects $V_-$. Such a subspace is
orthogonal to a one-dimensional subspace contained in $V_+$: we call \textit{polar
vector} of the complex line any generator of this subspace.  
\end{defi}

Note that a complex line is an isometric embedding of the complex hyperbolic
line $\HuC$. To any complex line $C$ is associated a unique holomorphic
involution fixing pointwise $C$, which we shall refer to as the
\textit{complex symmetry} with respect to $C$. The group PU(2,1) acts
transitively on the set of complex lines of $\HdC$.

\begin{defi}
We call \textit{flag} a pair $(C,p)$ where $C$ is a complex line and
$p$ is a point in $C\cap \partial \HdC$.
\end{defi}

\begin{lem}\label{transiflag}
PU(2,1) acts transitively on the set of flags of $\HdC$.
\end{lem}

\subsubsection{$\R$-planes}
\begin{defi}
An $\R$-plane is the intersection with $\HdC$ of the projection of a vectorial
Lagrangian subspace of $\C^{2,1}$.
\end{defi}

Every $\R$-plane $P$ is fixed
pointwise by a unique antiholomorphic isometric involution $I_P$,
which is the projectivization of the Langangian symmetry with respect
to any lift of $P$ as a vectorial Lagrangian. We will refer to $I_P$
as the \textit{Lagrangian reflection about $P$}. 
The standard example is the set of points of $\HdC$ with real
coordinates, which is fixed by the complex conjugation. We will refer
to this $\R$-plane as $\HdR\subset\HdC$. It is an embedding of the real
hyperbolic plane into $\HdC$.

As a consequence, we obtain
\begin{prop}
Let $Q$ be an $\R$-plane. There exists a matrix $M_Q\in$ SU(2,1) such
that 
\begin{equation}\label{relevR}M_Q\overline{M_Q}=1\mbox{, and }
I_Q\left(m\right)={\bf P}\left(M_Q\cdot \overline{{\bf m}}\right)\mbox{ for any
}m\mbox{ in }\HdC\mbox{ with lift }{\bf m}.\end{equation}  
\end{prop}
\begin{proof}
Let ${\bf Q}$ be a vectorial lift of $Q$, and choose $\R^3\subset
\C^{2,1}$ as a vectorial lift of $\HdR$. Since the group U(2,1) acts
transitively on the Lagrangian Grassmanian of $\C^{2,1}$, there
exists a matrix $A\in$ U(2,1) such that $A \R^3={\bf Q}$. The matrix
$M_Q=A\overline{A}^{-1}$ belongs to SU(2,1) and satisfies the condition
(\ref{relevR}). 
\end{proof}

\begin{rem}
If $I_1$ and $I_2$ are two Lagrangian reflections with associated
matrices ${\bf M}_1$ and ${\bf M}_2$, then their composition, which is
a holomorphic isometry, admits the matrix ${\bf
  M}_1\overline{{\bf M}}_2$ as a lift to SU(2,1).
\end{rem}
\subsection{Classical invariants}

\subsubsection{Invariant of two complex lines\label{2lines}}
\begin{defi}
Let $C_1$ and $C_2$ be two complex lines of $\HdC$, with polar vectors $\bc_1$ and $\bc_2$.
We set 
$$\phi(C_1,C_2)=\frac{|\la \bc_1,\bc_2\ra|^2}{\la \bc_1,\bc_1\ra \la \bc_2,\bc_2\ra}.$$
\end{defi}

\n Clearly, $\phi(C_1,C_2)$ does not depend on the choice of the lift
in the pair of polar vectors, and is PU(2,1)-invariant. We recall the geometric
interpretation of $\phi$, and we refer to \cite{Go} for details: 

\begin{itemize}
\item[-] $\phi(C_1,C_2)>1$ if  $C_1$ and $C_2$ are disjoint in
  $\HdC$. In this case, the distance $d$ between $C_1$ and $C_2$ is
  given by the formula $\phi(C_1,C_2)=\cosh^2(d/2)$. 
\item[-] $\phi(C_1,C_2)=1$ if $C_1$ and $C_2$ are either identical or
  asymptotic, by which we mean that they meet in $\partial \HdC$. 
\item[-] $\phi(C_1,C_2)<1$ if $C_1$ and $C_2$ intersect. The angle
  $\theta$ of their intersection is given by the relation
  $\phi(C_1,C_2)=\cos^2(\theta)$. 
\end{itemize} 

\n Note that two complex lines are orthogonal if and only if
$\phi(C_1,C_2)=0$. The $\phi$-invariant classifies 
pairs of distinct complex lines up to isometries.
\begin{prop}
Let $C_1,C_2,D_1,D_2$ be 4 complex lines such that $C_1\ne C_2$ and
$D_1\ne D_2$. There exists an isometry $g\in PU(2,1)$ such that $D_1=g
C_1$ and $D_2=g C_2$ if and only if $\phi(C_1,C_2)=\phi(D_1,D_2)$.  
\end{prop}
\begin{proof}
It is clear that if the two pairs are isometric, their invariant is
the same. Reciprocally, choose polar vectors $\bc_1,\bc_2,\bd_1,\bd_2$
with norm 1 and such that $\la \bc_1,\bc_2\ra $ and $\la
\bd_1,\bd_2\ra $ are in $\R_{\ge 0}$. 

 As $C_1\ne C_2$, the vectors $\bc_1$ and $\bc_2$ are independent and
 the same is true for $\bd_1$ and $\bd_2$. With the assumption on the
 $\phi$-invariant, the Gram matrices of $(\bc_1,\bc_2)$ and
 $(\bd_1,\bd_2)$ are identical. It means that one can find an isometry
 which maps $\bc_1$ on $\bd_1$ and $\bc_2$ on $\bd_2$. This ends the
 proof. 
\end{proof}

We will need the  following lemma.

\begin{lem}\label{stabilisateur}
Let $C_1$ and $C_2$ be two non orthogonal distinct complex lines,
and $p_1$ a point in $\partial\HdC\cap C_1$ which is not in $C_2$. Except for the identity, no isometry
preserves $C_1$ and $C_2$ and fixes $p_1$.
\end{lem}

\begin{proof}
Pick $\bc_1$ and $\bc_2$ two vectors polar to $C_1$ and $C_2$ of norm
1 and $\bp_1$ a lift of $p_1$. Writing $\la \bp_1,\bc_2\ra=a$  and $\la
\bc_1, \bc_2\ra=b$, the hermitian form has the following matrix in the
basis $\left(\bp_1,\bc_1,\bc_2\right)$
$$H=
\begin{bmatrix}
0 & 0 & a \\
0 & 1 & b \\
\bar a & \bar b & 1
\end{bmatrix}
$$
An isometry having the requested property has a diagonal lift to
SU(2,1) in this basis. The result is obtained by writing the isometry  
condition $^t\overline{M} H M=H$, and by using the fact that $b$ is
non zero since $C_1$ and $C_2$ are not orthogonal.
\end{proof}

\subsubsection{Invariants of three complex lines\label{3lines}}

Let $C_1,C_2$ and $C_3$ be three complex lines in $\HdC$. We will say that
they are in {\it generic position} if their polar vectors form a basis
of $\C^3$. There are three invariants of the triple $(C_1,C_2,C_3)$
given by the $\phi$-invariant of all pairs of complex lines. We will
need a fourth one (given in the following definition) to classify
all triples up to isometry.

\begin{defi}
Let $C_1,C_2,C_3$ be three complex lines in $\HdC$ with respective
polar vectors $\bc_1,\bc_2,\bc_3$. Then we set  
$$\Phi(C_1,C_2,C_3)=\frac{\la\bc_1,\bc_2\ra\la\bc_2,\bc_3\ra\la\bc_3,\bc_1\ra}
{\la\bc_1,\bc_1\ra\la\bc_2,\bc_2\ra\la\bc_3,\bc_3\ra}.$$
\end{defi}

\n The importance of this invariant should be clear from the following
two propositions (see \cite{Pra}):

\begin{prop}\label{troisdroitesrela}
Let $C_1,C_2,C_3$ be three complex lines of $\HdC$ in generic
position. For simplicity, we denote by $\phi_{ij}$ the
$\phi$-invariant of $C_i$ and $C_j$ and by $\Phi_{ijk}$ the
$\Phi$-invariant of $C_i,C_j,C_k$. These invariants enjoy the
following properties.
\begin{enumerate}
\item
For all distinct $i,j,k\in\{1,2,3\}$, the following relations are satisfied.
\begin{equation}\label{syminv}\phi_{ij}=\phi_{ji},
  \Phi_{ijk}=\Phi_{jki}=\ba{\Phi_{ikj}} \mbox{ and }
  \Phi_{ijk}\Phi_{ikj}=\phi_{ij}\phi_{jk}\phi_{ki}.
\end{equation}

\item The four invariants satisfy to the inequality 

 \begin{equation}\label{detgram}
   1-\phi_{12}-\phi_{23}-\phi_{31}+\Phi_{123}+\Phi_{132}<0.\end{equation} 
\end{enumerate}
\end{prop}
\begin{proof}
\begin{enumerate}
\item These relations are straightforward from the definitions of the
  invariants $\phi$ and $\Phi$.
\item Let $C_1,C_2,C_3$ be three complex lines in generic position. Let
$\bc_1,\bc_2,\bc_3$ be three polar vectors associated to these
lines. Let $G$ be the Gram matrix of the basis
$(\bc_1,\bc_2,\bc_3)$. A direct computation shows that the left-hand
side of relation (\ref{detgram}) is equal to 
$$\Delta(C_1,C_2,C_3)=\frac{\det G}{\la\bc_1,\bc_1\ra
  \la\bc_2,\bc_2\ra \la\bc_3,\bc_3\ra}.$$ 
This number is an invariant of the triple $(C_1,C_2,C_3)$. As the Gram
matrix represents the hermitian form, it has signature $(2,1)$
and its determinant is negative.
\end{enumerate}
\end{proof}

\begin{prop}\label{troisdroitesclass}
Consider non negative real numbers $\phi_{ij}$ and complex numbers $\Phi_{ijk}$
satisfying the relations of proposition \ref{troisdroitesrela}. There
exists a triple $C_1,C_2,C_3$ in generic position, unique up to
isometry, such that for all distinct  $i,j,k$ in $\{1,2,3\}$ the relations 
$\phi(C_i,C_j)=\phi_{ij}$ and $\Phi_{ijk}=\phi(C_i,C_j,C_k)$ hold.
\end{prop}  

\begin{proof}
  Consider real numbers $\phi_{ij}$ and complex numbers
$\Phi_{ijk}$ satisfying the relations (\ref{syminv}) and
(\ref{detgram}), and $\C^3$ with its canonical basis
$(e_1,e_2,e_3)$. We define a hermitian form $h$ on it by setting

$$\begin{array}{cc}h( e_i,e_i)=1\mbox{ for all $i$, } & h(e_1,e_2)
  =\sqrt{\phi_{12}},\\
\\
h(e_2,e_3)=\sqrt{\phi_{23}},
&
h(e_3,e_1)= \dfrac{\Phi_{123}}{\sqrt{\phi_{12}\phi_{23}}}.
\end{array}$$

The matrix of $h$ in the basis $(e_1,e_2,e_3)$ has unit diagonal
entries -- thus positive trace-- and according to the relation
(\ref{detgram}), it has negative determinant. As a consequence, $h$ has
signature $(2,1)$. By the classification of hermitian forms, this
model is conjugate to the standard one, and the vectors $e_1,e_2,e_3$
map to the polar vectors of the desired complex lines. Moreover, the
few choices we made disappear projectively, hence the triple of
complex lines is unique up to isometry. 
\end{proof}

\section{Invariants of  flags\label{flaginv}}
\subsection{Invariant of two flags}
\begin{defi}
Let $(C_1,p_1)$ and $(C_2,p_2)$ be two flags in $\HdC$. We will say
that they are in {\it generic position} if $p_1$ does not belong to
$C_2$, $p_2$ does not belong to $C_1$ and $C_1$ is not orthogonal to $C_2$.
\end{defi}

\begin{rem}
The condition of non-orthogonality of $C_1$ and $C_2$ will be needed to define the elementary isometries
associated to a triple of flags in a unique way (see propositions \ref{isomelem} and \ref{elemRplan}).
\end{rem}
\begin{defi}\label{defm}
Let $(C_1,p_1)$ and $(C_2,p_2)$ be two flags in generic position. Let
$\bc_1, \bc_2$ be polar vectors of $C_1,C_2$ and $\bp_1,\bp_2$ be
representatives of $p_1,p_2$. We set  
$$m[(C_1,p_1),(C_2,p_2)]=\frac{\la\bc_1,\bc_2\ra \la
  \bp_1,\bp_2\ra}{\la \bc_1,\bp_2\ra \la \bp_1, \bc_2 \ra}.$$ 
\end{defi}

 \n This invariant is a complex generalization of the $\phi$-invariant
of two complex lines. Its properties are summed up in the following
proposition. 

\begin{prop}
Let $(C_1,p_1)$ and $(C_2,p_2)$ be two flags in generic position, and
$m_{12}$ their invariant $m[(C_1,p_1),(C_2,p_2)]$. 
\begin{enumerate}
\item The two invariants $\phi(C_1,C_2)$ and $m_{12}$ are linked by the relation 
\begin{equation}\label{mphi}\phi(C_1,C_2)=\left|\frac{m_{12}}{m_{12}-1}\right|^2. 
\end{equation}    
\item For any complex number $m_{12} \in \C\setminus\{0,1\}$ there
exists a pair of flags $(C_1,p_1),(C_2,p_2)$ in generic position such that
$m[(C_1,p_1),(C_2,p_2)]=m_{12}$. This pair is unique up to isometry.  
\end{enumerate}
\end{prop}

\begin{proof}
\begin{enumerate}
\item
Let $(C_1,p_1)$ and $(C_2,p_2)$ be two flags in generic position,
$\bc_1, \bc_2$ be polar vectors of $C_1,C_2$ and $\bp_1,\bp_2$ be
representatives of $p_1,p_2$. The family of vectors
$(\bc_1,\bc_2,\bp_1,\bp_2)$ is linearly dependent, hence the determinant
of its Gram matrix vanishes. Computing this determinant and dividing
by non vanishing factors, we obtain the relation $|m_{12}-1|^2
\phi_{12}=|m_{12}|^2$. From that relation we see that $m_{12}$ cannot
be equal to 1 since the two flags are in generic position and therefore $\phi_{12}$ is non-zero. 
This proves relation (\ref{mphi}).
\item
In order to prove the last part of the proposition, we make the
following observation: given two flags $(C_1,p_1)$ and $(C_2,p_2)$ in
generic position, there is a unique complex line $C_3$ joining $p_1$
and $p_2$. Following proposition \ref{troisdroitesclass}, the triple
$(C_1,C_2,C_3)$ is determined by its $\phi$-invariants, hence,  we can
classify couples of flags using $\phi$-invariants. 
 
 More precisely, as $C_1$ and $C_3$ are asymptotic (they meet on
$p_1\in\partial\HdC$), their $\phi$ invariant $\phi_{13}$ equals
$1$. For the same reason, $\phi_{23}=1$. As a consequence of
relation (\ref{syminv}), we obtain the equality
$|\Phi_{123}|^2=\phi_{12}$. Plugging these values into the relation
(\ref{detgram}) yields   

\begin{equation}\Delta_{123}=-1-\phi_{12}+\Phi_{123}+\Phi_{132}=-|1-\Phi_{123}|^2.
\end{equation} 

Suppose that we have $\Phi_{123}\ne 1$, then the complex lines
$C_1,C_2,C_3$ are in generic position. Let $\bc_1,\bc_2,\bc_3$ be
polar vectors of these lines. They form a basis of $\C^3$ and the
linear forms $\la \cdot,\bc_1\ra,\la \cdot,\bc_2\ra,\la
\cdot,\bc_3\ra$ form a linear basis of the dual of $\C^3$. One can
find a unique anti-dual basis $\bd_1,\bd_2,\bd_3$ such that for all
$i,j$ in $\{1,2,3\}$ one has $\la \bd_i,\bc_j \ra=\delta_{ij}$. A
direct computation shows that the Gram matrix of the hermitian form in
the basis $(\bd_1,\bd_2,\bd_3)$ is the inverse of the Gram matrix of
$(\bc_1,\bc_2,\bc_3)$. Moreover $\bd_1$ being orthogonal to $\bc_2$
and $\bc_3$, it is a representative of $p_2$ and $\bd_2$ is a
representative of $p_1$. Using these  representatives we get (see remark \ref{hermcross} above)  
\begin{eqnarray*}
m_{12}&=&\la \bc_1,\bc_2\ra \la \bd_1,\bd_2\ra=\frac{\la \bc_1,\bc_2\ra
\left( 
\la \bc_3,\bc_1\ra\la \bc_2,\bc_3\ra-\la \bc_2,\bc_1\ra\la \bc_3,\bc_3\ra
\right)
}{\la \bc_1,\bc_1\ra\la \bc_2,\bc_2\ra\la
  \bc_3,\bc_3\ra\Delta_{123}}=\frac{\Phi_{123}-\phi_{12}}{\Delta_{123}}\\ 
&=& \frac{\Phi_{123}-|\Phi_{123}|^2}{-|1-\Phi_{123}|^2}=\frac{\Phi_{123}}{\Phi_{123}-1}.
\end{eqnarray*} 

This proves that $m_{12}=1$ if and only if $\Phi_{123}=1$ and that
$m_{12}$ classifies couples of flags as $\Phi_{123}$ does. 
\end{enumerate}
\end{proof}
\begin{rem}\label{hermcross}
Note that this anti-dual basis is usualy used in the literature under a slightly 
different form, using the so-called hermitian cross-product. The vector $\bd_2$ is proportionnal 
to the hermitian cross-product of $\bc_1$ and $\bc_3$, denoted by $\bc_1\boxtimes \bc_3$. It is a simple
computation using hermitian cross-product to check that $\la \bd_1,\bd_2\ra$ equals 
$\left(\la\bc_3,\bc_1\ra\la\bc_2,\bc_3\ra-\la\bc_2,\bc_1\ra\la\bc_3,\bc_3\ra\right)\cdot\left(\la\bc_1,\bc_1\ra\la\bc_2,
\bc_2\ra\la\bc_3,\bc_3\ra\Delta_{123}\right)^{-1}$.
See \cite{Go} for details.
\end{rem}

\subsection{Invariant of a flag and two complex lines}

\begin{defi} \label{defdelta}
Let $(C_1,p_1)$ be a flag and $C_2,C_3$ be two complex lines such that
the three complex lines $C_1,C_2$ and $C_3$ are in generic position
and such that $p_1$ does not belong to $C_2$ nor $C_3$. Take
$\bc_1,\bc_2,\bc_3$ three polar vectors of $C_1,C_2,C_3$ and $\bp_1$ a
representative of $p_1$. Then we set: 

\begin{equation}\delta[(C_1,p_1),C_2,C_3]=\delta^1_{23}=\frac{\la\bc_2,\bc_3\ra\la
    \bp_1,\bc_2\ra}{\la\bc_2,\bc_2\ra
    \la\bp_1,\bc_3\ra}.\end{equation} 
\end{defi}

\n This invariant may be viewed as a coordinate of $p_1$ knowing $C_1$,
$C_2$ and $C_3$. Its main properties are summed up in the following
proposition. 

\begin{prop}\label{1flag2droites}
Let $(C_1,p_1)$ be a flag and $C_2,C_3$ be two complex lines such that
the three complex lines $C_1,C_2$ and $C_3$ are in generic position
and such that $p_1$ does not belong to $C_2$ nor $C_3$. The invariants
$\delta^1_{23}$ and $\delta^1_{32}$ satisfy the following equations: 

\begin{eqnarray}
\phi_{23}& = &\delta^1_{23}\delta^1_{32}\label{deltadelta}\\
\nonumber\\
0 & = & (1-\phi_{13})|\delta^1_{23}|^2+2\Re\left[
  (\Phi_{132}-\phi_{23})\delta^1_{23} \right] +
\phi_{23}(1-\phi_{12})\label{deltacercle} 
\end{eqnarray}
\n Reciprocally, take $C_1,C_2,C_3$ three complex lines in generic
position. Any non zero value of $\delta^1_{23}$ which satisfies the second equation corresponds to
a unique point $p_1$ in $C_1$ which is not on $C_2$ nor on $C_3$. 
\end{prop}

\begin{proof}
The first equation is a direct consequence of the definition. For the second
one, let $(\bc_1, \bc_2, \bc_3)$ be  a basis of $\C^3$ formed by polar vectors for $C_1,C_2,C_3$.
 Let $(\bd_1,\bd_2,\bd_3)$ be its anti-dual basis. We will use the latter basis to prove relation
 (\ref{deltacercle}). We recall that the
matrix of the Hermitian form in the basis $(\bd_1,\bd_2,\bd_3)$ is the
inverse of the Gram matrix of $(\bc_1,\bc_2,\bc_3)$.  

As $p_1$ belongs to $C_1$, its representative is a linear combination
of $\bd_2$ and $\bd_3$, and we may thus write $\bp_1=a \bd_2+b \bd_3$. The
coordinates $a$ and $b$ can be recovered by computing the hermitian products $\la
\bp_1,\bc_2\ra=a$ and $\la \bp_1, \bc_3\ra=b$. In particular, this implies 
\begin{equation}\delta^1_{23}=\frac{\la \bc_2,\bc_3\ra a}{\la \bc_2,\bc_2 \ra b}.\label{projcoord}\end{equation} 

By expressing that $\bp_1$ is in the isotropic cone of the Hermitian
form,  we obtain the relation (\ref{deltacercle}).   

On the other hand, if we know $\delta^1_{23}$, then according to relation (\ref{projcoord}),
 we know projective coordinates for $\bp_1$. If $\delta^1_{23}$ satisfies (\ref{deltacercle}), the 
 vector $\bp_1$ must be on the cone of the quadratic
form. It proves that $\delta^1_{23}$ determines the position of $p_1$
on $C_1$ as asserted. 
\end{proof}

\subsection{Summary : invariants of three flags}

In the remaining part of the article, we will be interested in the space
of configurations of three flags. Let us sum up what are the relevant
invariants for such configurations. 

\begin{defi}\label{genflag}
We will say that three flags $\left(C_i,p_i\right)_{i=1,2,3}$ are in
\textit{generic position}  if they are pairwise in generic position,
and if the triple of complex lines $(C_1,C_2,C_3)$ is also in generic
position, that is, if
\begin{itemize}
\item any two of the complex lines are disctinct and non-orthogonal,
\item any triple of vectors polar to the $C_i$'s is a basis of $\C^3$.
\end{itemize}
\end{defi}

\n We classify now the triples of flags up to PU(2,1).

\begin{theo}\label{classdrap}
Let $(C_1,p_1)$, $(C_2,p_2)$ and $(C_3,p_3)$ be three flags in
generic position.  The configuration of these flags modulo
holomorphic isometry is classified by the invariants
$\phi_{ij},\Phi_{ijk}$ and $\delta^i_{jk}$ for 
all distinct $i,j,k$ in $\{1,2,3\}$. 
These invariants satisfy the following equations for all $i,j,k$:
\begin{itemize}
\item[(\ref{syminv})]
  \begin{center}$\quad\phi_{ij}=\phi_{ji}=\ba{\phi_{ij}}>0$,
    $\Phi_{ijk}=\Phi_{jki}=\ba{\Phi_{ikj}}$ and
    $\Phi_{ijk}\Phi_{ikj}=\phi_{ij}\phi_{jk}\phi_{ki}$.\end{center} 
\item[(\ref{detgram})]
  \begin{center}$\quad\Delta_{ijk}=1-\phi_{ij}-\phi_{jk}-\phi_{ki}+\Phi_{ijk}+\Phi_{ikj}<0.$\end{center}  
\item[(\ref{deltadelta})]\begin{center}
    $\quad\delta^i_{jk}\delta^i_{kj}=\phi_{ij}.$\end{center} 
\item[(\ref{deltacercle})]\begin{center}
    $\quad(1-\phi_{ik})|\delta^i_{jk}|^2+2\Re\left[
      (\Phi_{ikj}-\phi_{jk})\delta^i_{jk} \right] +
    \phi_{jk}(1-\phi_{ij})=0.$\end{center} 
\end{itemize}
The space of solutions is a manifold of dimension 7. Moreover, the
invariants $m_{ij}$ attached to pairs of flags are expressed in terms
of the other invariants as follows :
\begin{eqnarray}
m_{ij}\Delta_{ijk}\phi_{ik}\phi_{jk} & = & \phi_{ik}\phi_{jk}(\Phi_{ijk}-\phi_{ij})+
\phi_{ik}(\phi_{ij}\phi_{jk}-\Phi_{ijk})\delta^i_{kj}\nonumber\\
& &+\phi_{jk}(\phi_{ij}\phi_{ik}-\Phi_{ijk})\ba{\delta^j_{ki}}+
\Phi_{ijk}(1-\phi_{ij})\delta^i_{kj}\ba{\delta^j_{ki}}\label{metautres}
\end{eqnarray}
\end{theo}

\begin{proof}
The first part of the proof is nothing but a summary of the preceding
sections. Let us now compute $m_{12}$. The two other $m$-invariants
are obtained in the same way.  Choose
$\bc_1,\bc_2,\bc_3$ polar vectors of $C_1,C_2,C_3$ and let
$(\bd_1,\bd_2,\bd_3)$ be the anti-dual basis as usual. Then, using the proof of proposition
\ref{1flag2droites}, one can find explicit  coordinates for representatives of $p_1$ and $p_2$ in
the basis $(\bd_1,\bd_2,\bd_3)$.  Precisely, we can choose 
$$
\left\{\begin{array}{cc}
\bp_1=\la \bc_3,\bc_2\ra \bd_2 +\la \bc_3,\bc_3\ra \delta^1_{32}\bd_3 \\
 \\
\bp_2=\la \bc_3,\bc_1\ra \bd_1 +\la  \bc_3,\bc_3\ra \delta^2_{31} \bd_3. \\
\end{array}\right.$$

To obtain a formula for $m_{12}$, we just need to replace $\bp_1$ and
$\bp_2$ in the definition of $m_{12}$ by the expressions above. 
We obtain the relation (\ref{metautres}) after a computation.
\end{proof}

\section{Elementary isometries associated to a triple of flags\label{stand}}
\subsection{$\R$-planes associated to a triple of flags and elementary isometries}
In this paragraph, we define the elementary isometries associated to a triple of
flags. More precisely, we prove the

\begin{prop}\label{isomelem}
Let $F_i=\left(C_i,p_i\right)$ for $i=1,2,3$ be a triple of flags in
generic position such that any two complex lines are not asymptotic.
\begin{enumerate}
\item For any pair $(i,j)$ with $i\neq j$, there exists a unique
  isometry $E_{ij}$ exchanging $C_i$ and $C_j$,
  and mapping $p_j$ to $p_i$. It is called the \textit{exchange isometry}
    associated to the pair of flags $F_i$ and $F_j$.
\item There exists a unique isometry $T^i_{jk}$ fixing $p_i$ and
  preserving $C_i$ which maps $C_k$ to a complex line $C'_k$
  satisfying $R_{C'_k}(p_i)=R_{C_j}(p_i)$, where $R_C$ is the complex symmetry with respect 
  to the complex line $C$.
   It is called the \textit{transfer isometry}
    associated to the ordered triple of flags $\left(F_i,F_j,F_k\right)$.
\end{enumerate}
\end{prop}

We will give a geometric proof of this proposition, showing that the exchange and transfer isometry are
obtained as products of Lagrangian reflections which are canonically associated to a triple of
flags satisfying the assumption of proposition \ref{isomelem}. 

\begin{prop}\label{elemRplan}
Let $C_1$ and $C_2$ be two complex lines which are neither orthogonal
nor asymptotic.
\begin{enumerate}
\item Let $p_1$ be a point in $\partial C_1$. There exists a unique
  $\R$-plane $P$ such that $I_P$, the inversion in $P$, preserves both
  $C_1$ and $C_2$, and fixes $p_1$.
\item Let $p_2$ be a point in $\partial C_2$. There exists a unique
  $\R$-plane $Q$ such that $I_Q$, the inversion in $Q$, swaps $C_1$
  and $C_2$ and maps $p_1$ to $p_2$.
\item Let $m$ and $n$ be two points in the boundary of $\HdC$, not
  belonging to $\partial C_1$. There exists a unique Lagrangian reflection
  preserving $C_1$ and swapping $m$ and $n$. 
\item Let $p_1$, $p_2$ and $p_3$ be three points of $\partial\HdC$, not contained in the
boundary of a complex line. There exists a unique Lagrangian reflection fixing $p_1$ and
swapping $p_2$ and $p_3$.
\end{enumerate} 
\end{prop}

\begin{proof}
Let ${\bf c}_k$ be a polar vector for $C_k$ normalized so that $\la{
  \bf c}_k,{\bf c}_k\ra=1$. Let $\bp_1$ be a lift of $p_1$.
  Rescaling if necessary, we may assume that both
  $a=\la\bp_1,\bc_2\ra$ and $b=\la\bc_1,\bc_2 \ra$ are real 
  (in fact $\la\bc_1,\bc_2\ra$ is equal to $\sqrt{\phi_{12}}$).
\begin{enumerate}
\item The hermitian form admits in the basis $\left({\bf p_1},{\bf c}_1,{\bf
    c}_2\right)$ the matrix 
$$H=
\begin{bmatrix}
0 & 0 & a\\
0 & 1 & b \\
a & b & 1 \\
\end{bmatrix}$$
The hermitian product $b$ is non-zero since $C_1$ and $C_2$ are non-orthogonal. 
In this basis, any lift of a Lagrangian reflection fixing $p_1$ and preserving
$C_1$ and $C_2$ must be diagonal. It follows after writing the isometry
condition $M^* H M=H$ that there is only one such reflection, given in
this basis by ${\bf m}\longrightarrow \overline{\bf m}$. 
\item 
This time, we use the basis $\left({\bf c}_1,{\bf c}_2,{\bf d}\right)$, where $\bf d$ a vector orthogonal to
$\bc_1$ and $\bc_2$ with norm $b^2-1$ (indeed, we are setting $\bf d={\bf c}_1\boxtimes {\bf c}_2$, see remark
\ref{hermcross}). 
The hermitian form is given by the matrix
$$H=\begin{bmatrix}
1 & b & 0\\
b & 1 & 0 \\
0 & 0 & b^2-1  \\
\end{bmatrix}\quad \mbox{($|b|=1$ iff $C_1$ and $C_2$ are asymptotic)}
$$

We may choose the lifts of $p_1$ and $p_2$ as follows :
$${\bf p}_1=\begin{bmatrix} - b \\ 1 \\
  e^{i\theta_1}\end{bmatrix}\mbox{ and }
{\bf p}_2=\begin{bmatrix} 1 \\ -b \\
  e^{i\theta_2}\end{bmatrix}\mbox{ with }\theta_i\in\R. 
$$
The fact that $I_Q$ exchanges $C_1$ and $C_2$ implies that any matrix for $I_Q$ has the
form 
$$
\begin{bmatrix}
0 & \alpha & 0 \\
\beta & 0 & 0\\
0 & 0 & \gamma \\
\end{bmatrix}.
$$
Writing the isometry condition and the fact that
$I_Q\left(p_1\right)=p_2$, provides relations determining $\alpha$,
$\beta$ and $\gamma$. The result
follows. 
\item We may choose lifts $\bm$ and $\bn$ of $m$ and $n$ such that $\la\bm,\bn\ra=1$ and a
unit vector $\bc$ polar to the complex line containing $m$ and $n$. In the basis
$(\bm,\bc,\bn)$, where the hermitian form has matrix $J$, the complex line $C_1$ is polar
to some vector $\bc_1=\begin{bmatrix} \alpha & \beta &\gamma \end{bmatrix}^T$. It is a direct
computation to check that a Lagragian reflection swapping $m$ and $n$ and preserving $C_1$
lifts to the matrix below. Hence, it exists and is unique.
$$
\begin{bmatrix}
0 & 0 & \alpha/\bar\gamma\\
0 & \beta/\bar\beta & 0 \\
\gamma/\bar\alpha & 0 & 0
\end{bmatrix}$$
\item In the proof of the previous item, we have not used the fact that $\bc_1$ was a
positive vector. Thus the same result as {\it 3} remains true if we change $C_1$ to a
boundary point, that is, $\bc_1$ to a null vector. If the three points are in a 
complex line, then we lose the uniqueness.
Note that this fourth part of the proposition is classical (see for instance lemma 7.17 of
\cite{Go})
\end{enumerate}
\end{proof}

\begin{proof}[Proof of proposition \ref{isomelem}]
\begin{enumerate}
\item  Let $h_1$ and $h_2$ be two isometries having the requested
  properties. Then $h_2^{-1}\circ h_1$ preserves both $C_i$ and $C_j$,
  and fixes $p_i$. According to the lemma \ref{stabilisateur}, this implies that
  $h_2$ and $h_1$ are equal. This proves the uniqueness. To prove the existence part, 
  we apply the first two items of proposition \ref{elemRplan}.
  \begin{itemize}
  \item There exists a unique Lagrangian reflection $I_2$ preserving $C_i$ and $C_j$ and fixing $p_i$
  (this follows from part 1 of proposition \ref{elemRplan}).
  \item There exists a unique Lagrangian reflection $I_1$ swapping $C_i$ and $C_j$, and exchanging $p_i$
  and $I_2(p_j)$. This is part 2 of proposition \ref{elemRplan}, which may be applied since $I_2(p_j)$
  belongs to $C_j$.
    \end{itemize}
  The isometry $E_{ij}=I_1\circ I_2$ has the requested properties.
  
  \item The uniqueness is proved in the same way as for \textit{1}.
  To prove the existence, we apply the third and fourth part of proposition
  \ref{elemRplan}.
  \begin{itemize}
  \item The two points $R_{C_3}(p_1)$ and $R_{C_2}(p_1)$ do not belong to $\partial C_1$ since the three
  complex lines are non-asymptotic. Thus, there exists a unique Lagrangian reflection 
  $I_3$ preserving $C_1$ and swapping $R_{C_3}(p_1)$ and $R_{C_2}(p_1)$ (this follows from part \textit{3} 
  of proposition \ref{elemRplan}). Note that $I_3$ does dot fix $p_1$.
  \item  The three points $p_1$, $I_3(p_1)$ and $R_{C_2}(p_1)$ do not belong to a common
  complex line, for else $C_1$ and $C_2$ would be asymptotic. Thus we may apply the fourth
  part of proposition \ref{elemRplan} to obtain a (unique) Lagrangian reflection $I_4$ fixing
  $R_{C_2}(p_1)$, and swapping  $p_1$ and $I_3(p_1)$.
    \end{itemize}
  The isometry $I_4\circ I_3$ has the requested properties (note that since $I_4$ swaps
  $p_1$ and $I_3(p_1)$ which both belong to $C_1$, it preserves $C_1$).
  \end{enumerate}
\end{proof}

\subsection{Standard position of a triple of flags and elementary isometries}

\begin{defi}
\begin{itemize}
\item[-] Let $(C_1,p_1)$ be a flag and $C_2$ be a complex line. We will say
that they are in {\it generic position} if $p_1$ does not belong to
$C_2$, and if $C_1$ and $C_2$ are distinct and non-orthogonal. 
\item[-]
We say that $(C_1,p_1)$ and $C_2$ are in \textit{standard position} if
$p_1,C_1$ and $C_2$  are respectively represented by the following
vectors: 
$$\bp_1=\mat{1\\0\\0}, \bc_1=\mat{0\\1\\0}, \bc_2=\mat{a\\ \sqrt{2}\\
  1}\text{ for }a\in (-1,+\infty).$$ 
\end{itemize}
\end{defi}

\n The condition on $C_2$ is equivalent to saying that $R_{C_2}(p_1)$ is represented by the vector 
$\mat{-1 & \sqrt{2} & 1}^T$.
The motivation for this definition is the following proposition:

\begin{prop}
Let $(C_1,p_1)$ be a flag and $C_2$ be a complex line in generic
position. 
\begin{itemize}
\item[-] There exists a unique couple in standard position which is
  isometric to $\left((C_1,p_1),C_2\right)$.

\item[-] The parameter $a$ is given by $\phi(C_1,C_2)=(1+a)^{-1}$
\end{itemize}
\end{prop}

\begin{proof}
Since $PU(2,1)$ acts transitively on the set of flags of $\HdC$, we can
assume that $\bp_1$ and $\bc_1$ are in standard position. The
isometries $g$ in PU(2,1) stabilizing the standard flag admit lifts
to SU(2,1) of the following form :  
$${\bf g}=\begin{bmatrix}\lambda & 0 & it\lambda\\ 0& \ba{\lambda}/\lambda& 0\\ 0 & 0 &
  1/\ba{\lambda}\end{bmatrix}\mbox{ with
}\lambda\in\C\setminus\{0\}\mbox{ and }t\in \R.$$  
Note that $\lambda$ is well-defined up to multiplication by a cubic
root of 1. Now, a generic polar vector for $C_2$ and its image by $g$
are given by

$$\bc_2=\mat{a\\b\\1}\mbox{ and
}{\bf g}\bc_2\sim\mat{|\lambda|^2(a+it)\\\ba{\lambda}^2b/\lambda\\1}\mbox{ with
}|b|^2+2\Re(a)>0.$$ 

The assumption that $\bc_1$ and $\bc_2$ are not orthogonal, implies that $b\ne 0$. 
This means that there is only one isometry which
stabilizes the standard flag and maps $\bc_2$ in standard
position. Namely, we have to set $t=-\Im(a)$ and solve
$\ba{\lambda}^2b=\sqrt{2}\lambda$. This equation has three solutions
in $\lambda$ which represent the same element in $PU(2,1)$. The value
of $\phi(C_1,C_2)$ is given by a straightforward computation.
\end{proof}

\begin{rem}
Given three flags $(C_1,p_1)$, $(C_2,p_2)$, $(C_3,p_3)$, we can decide
to put $(C_1,p_1)$ and $C_2$ in standard position. However, we could have
chosen $(C_1,p_1)$ and $C_3$ or $(C_2,p_2)$ and $C_1$. All these
configurations can be obtained one from the other by applying
elementary isometries to the configuration.

As an example, assume that $(C_1,p_1)$ and $C_2$ are in standard
position, and apply the exchange $E_{12}$ isometry swapping $C_1$ and
$C_2$ and mapping $p_2$ to $p_1$. Their images
$\left(E_{12}\left(C_2\right),E_{12}\left(p_2\right)\right)$ and 
$E_{12}\left(C_1\right)$ are in standard position.

In the same way, applying the transfer isometry $T^1_{23}$ to the
triple $(C_1,p_1)$, $(C_2,p_2)$ and $(C_3,p_3)$ with $(C_1,p_1)$ and
$C_2$ in standard position makes $(C_1,p_1)$ and $C_3$ in standard
position. 
\end{rem}

\begin{prop}\label{isomflag}
Let $(C_1,p_1)$, $(C_2,p_2)$, $(C_3,p_3)$ be a triple of flags in
generic position and $\Theta:\C\to\C$ be the map defined by $\Theta(\rho
e^{i\theta})=\rho e^{i\theta/3}$ for $\rho\in [0,+\infty)$ and
$\theta\in (-\pi,\pi]$. 
 Assume that $(C_1,p_1)$ and $C_2$ are in standard position.
\begin{enumerate}
\item   
The transfer isometry $T^1_{23}$ is given by its lift to SU(2,1):
$${\bf T}^1_{23}=\mat{\mu & 0 & it\mu\\ 0&
  \ba{\mu}/\mu& 0\\ 0 & 0 & 1/\ba{\mu}}\mbox{ where
}\mu=\Theta(\frac{\delta^1_{23}\phi_{13}}{\Phi_{123}})\mbox{ and
}t=\Im\left(\frac{2\delta^1_{23}(\phi_{23}-\Phi_{132})}{\phi_{12}\phi_{23}}\right)
$$  
\item The exchange isometry $E_{12}$ is given by its lift to SU(2,1):
$$\bf{E}_{12}=\mat{
\dfrac{\lambda(z-\ba{z}-|z|^2)}{4|z(z-1)|^2} &
\dfrac{\sqrt{2}\ba{z}\lambda(z-\ba{z}-|z|^2)}{4|z(z-1)|^2}
+\dfrac{\lambda}{\sqrt{2}(z-1)}&
\dfrac{\lambda}{1-z}+\dfrac{\lambda(z-\ba{z}-|z|^2)^2}{4|z(z-1)|^2}\\ 
\dfrac{\ba{\lambda}}{\sqrt{2}\lambda(\ba{z}-1)} &
\dfrac{\ba{\lambda}}{\lambda (\ba{z}-1)} &
\dfrac{\ba{\lambda}(|z|^2-z-\ba{z})}{\lambda(\ba{z}-1)\sqrt(2)}\\ 
\dfrac{1}{\ba{\lambda}}& \dfrac{\sqrt{2}\ba{z}}{\ba{\lambda}}&\dfrac{-|z|^2+z-\ba{z}}{\ba{\lambda}}
}$$
where $z=1/\ba{m_{12}}$ and $\lambda=2\Theta(z(z-1))$. 

\end{enumerate}

\end{prop}

\begin{proof}
\begin{enumerate}
\item 
Suppose that the triple of flags $(C_1,p_1)$, $(C_2,p_2)$ and
$(C_3,p_3)$ is in generic position as it is specified in the
proposition, and suppose moreover that $(C_1,p_1)$ and $C_2$ are in
standard position. We can choose polar vectors $\bc_1,\bc_2,\bc_3$ and
representatives $\bp_1,\bp_2,\bp_3$ such that 
$$\bp_1=\mat{1\\0\\0},\,
\bc_1=\mat{0\\1\\0},\, \bc_2=\mat{1/\phi_{12}-1\\ \sqrt{2}\\ 1}.$$  
The
matrix we are interested in stabilizes $C_1$ and $p_1$ and sends $C_3$
to a standard complex line with polar vector
$$\bc_3'=\mat{1/\phi_{13}-1\\ \sqrt{2}\\ 1}.$$ 
Call ${\bf g}$ the inverse of the expected matrix, and compute the image of
$\bc'_3$ by ${\bf g}$:
$${\bf g}=\mat{\lambda & 0 & it\lambda\\ 0& \ba{\lambda}/\lambda& 0\\ 0 &
  0 & 1/\ba{\lambda}} \mbox{ and }\bc_3={\bf g}\bc_3'=\mat{\lambda(1/\phi_{13}-1+it)\\
  \ba{\lambda}\sqrt{2}/\lambda \\ 1/\ba{\lambda}}.$$
Computing explicit expressions for $\delta^1_{23},\phi_{23}$ and
$\Phi_{123}$ yields  equations for $\lambda$ and $t$. A direct
computation gives the formulas of the proposition. 

\item
The second matrix is obtained in three steps: let $(C_1,p_1)$ and
$(C_2,p_2)$ be two flags in generic position such that $(C_1,p_1)$ and
$C_2$ are in standard position. We look for a transformation which
sends $(C_2,p_2)$ and $C_1$ to a standard position. We find
explicitely a first transformation which sends $p_2$ to $p_1$. Then we
compose it with a Heisenberg translation (see remark \ref{heis} below) which sends the image of
$C_2$ by the first transformation to $C_1$. It remains to find a
matrix as in the first part which stabilize the standard flag
$(C_1,p_1)$ and sends the image of $C_1$ by the two first
transfomations to a standard complex line. The composition of these
matrices gives the formula of the proposition. 
\end{enumerate} 
\end{proof}
\begin{rem}\label{heis}
A Heisenberg translation is a unipotent parabolic isometry, given by the matrix
$$
\begin{bmatrix}
1 & -\bar w \sqrt{2}& -|w|^2+i\tau\\
0 &1 & w\sqrt{2}\\ 0 & 0 & 1
\end{bmatrix}\mbox{ with }w\in\C \mbox{ and }\tau\in\R.
$$
It is an element of the maximal unipotent subgroup of PU(2,1) fixing the vector $\begin{bmatrix}1 & 0 &
0\end{bmatrix}^T$,
which is a copy of the Heisenberg group of dimension 3.
\end{rem}

\begin{rem}
If $\left(C_1,p_1\right)$ and $p_2$ are in standard position, then the
$\R$-plane provided by the first part of proposition \ref{elemRplan}
is $\HdR$. The inversion in that plane is associated to the identity
matrix. As a consequence of proposition \ref{isomelem}, the associated
exchange isometry admits a lift of the form $M_1\circ
\overline{Id}=M_1$ where $M_1$ is the matrix of a Lagrangian
reflection. This shows that ${\bf E}_{12}\overline{{\bf E}_{12}}=1$.
\end{rem}

\section{Decorated triangulations and representations of $\pi_{g,p}$\label{deco}}
In this section, we will prove the theorem \ref{bij} that is stated in the introduction.\\

We denote by $\pi_{g,p}$ be the fundamental group of $\Sigma_{g,p}$, a surface of genus $g$
with $p$ punctures $x_1,\ldots x_p$, assuming $p>0$. Recall that $\widehat{\Sigma}_{g,p}$ is the universal covering of $\Sigma_{g,p}$. 
Provided that the inequality $2-2g-p<0$ is satisfied, the surface $\widehat{\Sigma}_{g,p}$ is homeomorphic to a topological disk 
and the punctures lift to a $\pi_{g,p}$-invariant subset $X$ of the boundary of $\widehat{\Sigma}_{g,p}$.

\begin{defi}
We set
 $$\mathfrak{R}_{g,p}=\left\{(\rho,F)\right\}/PU(2,1),$$
where $\rho$ is a morphism from $\pi_{g,p}$ to PU(2,1) and $F$ is a map from $X$ to the set of flags in $\HdC$ such that 
for any $x\in X$ and $g\in\pi_{g,p}$ one has $F(g.x)=\rho(g).F(x)$.
The group $PU(2,1)$ acts on $F$ by isometry on the target and acts on $\rho$ by conjugation: this action corresponds to
changing the base point in $\pi_{g,p}$. 
\end{defi}

For convenience, let us recall what will be called a triangulation of
$\Sigma$, which is sometimes referred to as an ideal triangulation. A
triangulation of $\Sigma$ is an oriented finite 2-dimensional
quasi-simplicial complex $T$ with an homeomorphism $h$ from the topological
realization $|T|$ of $T$ to $\Sigma$ which maps vertices to
punctures. By quasi-simplicial, we mean that two distinct triangles of $T$ can share the same vertices. 
By a slight abuse of notation, we will nevertheless refer to a 2-simplex by its vertices.

Given a triangulation $T$ of $\Sigma$, we can lift it to a triangulation of $\widehat{\Sigma}_{g,p}$. We thus obtain a triangulation of a disk with vertices on the boundary. 
Such a triangulation is isomorphic to the Farey triangulation which is a very nice and visual object (see \cite{FockGon2}). We may think that any triangulated surface is a quotient of the Farey triangulation.
Given a pair $(\rho,F)$ and a triangle $\Delta$ of $T$, we can pick a lift of $\Delta$ which has three vertices $x,y$ and $z$ in $X$. 

\begin{defi}
We will say that the pair $(\rho,F)$ is {\it generic with respect to $T$} if for any 
lifts of triangles of $T$ with vertices $x,y$ and $z$, the triple of flags 
$(F(x),F(y),F(z))$ is generic in the sense of definition \ref{genflag}. 
We denote by $\mathfrak{R}_{g,p}^T$ the subset of $\mathfrak{R}_{g,p}$ made of pairs 
which are generic with respect to $T$.
\end{defi}

\begin{defi}\label{defdeco}
Let $T$ be a triangulation of $\Sigma$. We denote by $\boX(T)$ the set
of triples $\left(\phi,\Phi,\delta\right)$ where : 
\begin{itemize}
\item $\phi$ is an $\R_{>0}$-valued function defined on the set of
  unoriented edges of $T$,
\item $\Phi$ and $\delta$ are $\C$-valued functions defined on the set of
  ordered faces of $T$.
\end{itemize}
From these data, we define auxiliary invariants in the following
way. For any ordered face $\left(i,j,k\right)$ of $T$, we set: 
\begin{eqnarray}
\Delta_{ijk}&=&1-\phi_{ij}-\phi_{jk}-\phi_{ik}+\Phi_{ijk}+\Phi_{ikj}\\
m^k_{ij}&=&\frac{1}{\Delta_{ijk}\phi_{ik}\phi_{jk} }\big[ \phi_{ik}\phi_{jk}(\Phi_{ijk}-\phi_{ij})+\nonumber\\
&&\phi_{ik}(\phi_{ij}\phi_{jk}-\Phi_{ijk})\delta^i_{kj}+\phi_{jk}(\phi_{ij}\phi_{ik}-\Phi_{ijk})
\ba{\delta^j_{ki}}+\Phi_{ijk}(1-\phi_{ij})\delta^i_{kj}\ba{\delta^j_{ki}}\big]   
\end{eqnarray}
The maps $\phi,\Phi$ and $\delta$ must satisfy the following relations for all ordered face $(i,j,k)$ in $T$: 
\begin{eqnarray}
\vert\Phi_{ijk}\vert^2 & =
&\phi_{ij}\phi_{jk}\phi_{ki}\label{modPhi }\\  
\Phi_{ijk} & = & \Phi_{jki}=\overline{\Phi_{ikj}} \label{sym3}\\
\Delta_{ijk}&<&0\\
0 & = & \vert\delta^{i}_{jk}\vert^2\left(1-\phi_{ik}\right)+
 2\Re\lbrack\delta^{i}_{jk}\left(\Phi_{ikj}-\phi_{jk}\right)\rbrack 
+\phi_{jk}\left(1-\phi_{ik}\right)
\label{cercledelta}  
\end{eqnarray}

\n Moreover, for any edge $(i,j)$ belonging to the faces $(i,j,k)$ and $(i,j,l)$, we impose the relation 
\begin{equation}m^k_{ij}=m^l_{ij}.\label{compatibilite}\end{equation}

\end{defi}

Before starting the proof of theorem \ref{bij}, let us give some useful constructions:

\begin{defi}\label{triangul}
Let $T$ be a triangulation of $\Sigma_{g,p}$. By definition, $T$ is a
quasi-simplicial 2-complex and there is a homeomorphism $h$ from $|T|$ to
$\Sigma_{g,p}$. We consider the following sub-complex of $|T|$:  
\begin{itemize}
\item vertices are combinations $V_{xy}=\frac{2}{3}x+\frac{1}{3}y$
  where $x$ and $y$ belong to the same edge in $T$,  
\item there are two types of simplicial edges: one from $V_{xy}$ to
  $V_{yx}$ for any edge $(x,y)$, and one from $V_{xy}$ to $V_{xz}$ for
  any two adjacent edges $(x,y)$ and $(x,z)$.
\item In each face of $T$, the edges constructed above draw an
  hexagon: we add to the sub-complex the corresponding 2-cell. 
  
  We denote by $HT$ and call {\it hexagonation} of $T$ the sub-complex we
  have obtained. It has the structure of a 2-dimensional CW-complex
  homeomorphic to $\Sigma_{g,p}$. 
\end{itemize}

\end{defi}
\begin{figure}
\begin{center}
\vspace{1cm}
\begin{pspicture}(-2,-2)(4,4)
\includegraphics[width=5cm,height=5cm,scale=0.4]{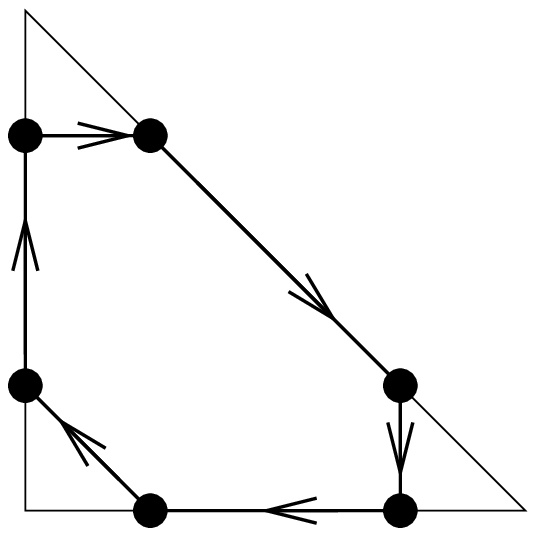}
\rput(0.2,0.1){$y$}
\rput(-5.2,0.1){$z$}
\rput(-5,5.2){$x$}
\rput(-2.5,-0.3){$E_{yz}$}
\rput(-2.3,3){$E_{xy}$}
\rput(-5.4,2.5){$E_{zx}$}
\rput(-4,3.2){$T^x_{zy}$}
\rput(-1.8,0.8){$T^y_{xz}$}
\rput(-3.6,0.8){$T^z_{yx}$}
\end{pspicture}
\vspace{-1cm}
\begin{center}
\caption[texte]{\protect\parbox[t]{12cm}{A hexagon associated to the triangle $(x,y,z)$, and
  elementary matrices associated to its sides.\label{simplichem}}}\end{center}
\end{center}
\end{figure}

 Let $T$ be a decorated triangulation of $\Sigma_{g,p}$. We will define from these data 
 a 1-cocycle $A$ in $Z^1(HT,$PU(2,1)). Let $s$ be an
oriented edge of $HT$. Associate to $s$ an
elementary matrix $A_s$ as follows:
\begin{itemize}
\item If $s=\left(V_{xy},V_{yx}\right)$ for some adjacent vertices of $T$, we set $A_s=E_{12}$
  where we replaced $m_{12}$ by $m_{xy}$. 
\item If $s=\left(V_{xy},V_{xz}\right)$, then we set $A_s=T^1_{23}$
  where we replaced all invariants by the decorations corresponding to
  the bijection $1\to x, 2\to y, 3\to z$. 
\end{itemize}

\begin{lem}\label{cocycle}
Let $T$ be a decorated triangulation of $\Sigma_{g,p}$. The mapping
$s\longrightarrow A_s$ is a 1-cocycle of $HT$ with values in PU(2,1).
\end{lem}
\begin{proof}
If $\left(x,y,z\right)$ is a face of $T$ and if $s_1\cdots s_6$ are
the sides of the associated hexagon, the product $\prod_{i=1}^6
A_{s_i}$ corresponds to an isometry of $\HdC$ stabilizing a flag and a
complex line. Hence, it is the identity map of $\HdC$ (see lemma \ref{stabilisateur}).
\end{proof}
We now go to the proof of theorem \ref{bij}.

\begin{proof}[Proof of theorem \ref{bij}]
We can finally prove the theorem by describing two mappings inverse one of each
other. For this purpose, fix a triangulation $T$ of $\Sigma_{g,p}$.

First, we associate to a decoration of $T$ a representation of
$\pi_{g,p}$ in PU(2,1) and an equivariant map $F$. Assume that $T$ is equipped with a decoration
$\left(\phi,\Phi,\delta\right)$ and choose a vertex $v=V_{a,b}$ of $HT$ as base point for the fundamental group of $\Sigma_{g,p}$.

Any loop $l$ of $\pi_1\left(\Sigma_{g,p},v\right)$ is homotopic to
a sequence $s=s_1,\ldots s_{k}$ of oriented edges of $HT$. One can
associate to $l$ the element of PU(2,1) corresponding to the product   
$$A_{s_{k}}\cdots A_{s_{1}}.$$

Because of the cocycle condition  given in lemme \ref{cocycle} above, this isometry does not depend
on the choice of the simplicial path homotopic to $l$. This gives rise to a
representation $\rho$ of $\pi_1(\Sigma_{g,p},v)$ into PU(2,1). Let us
now construct the map $F$. 
The choice of base point $v=V_{a,b}$ gives naturally a preferred lift of $a$ and $b$ in 
$\widehat{\Sigma}_{g,p}$ that we denote by 
$\widehat{a}$ and $\widehat{b}$ respectively. We choose $F(\widehat{a})$ and $F(\widehat{b})$ such that 
they are in standard position. Next, any element $x$ of $X$ is parametrized by a path from $v$ to a vertex of $HT$. We can suppose that this path $\gamma$ is simplicial. 
In that way, we set $F(x)=A_{\gamma}^{-1}.F_0$ where $F_0$ is the standard flag given by the vectors 
$\bc_1=\mat{0 & 1 &0}^T, \bp_1=\mat{1 &0 &0}^T$. 
One checks easily that this map $F$ is equivariant and generic with respect to $T$ and hence, the couple $(\rho,F)$ gives an element of $\mathfrak{R}^T_{g,p}$.

Conversely, given a couple $(\rho,F)$ generic with respect to $T$, we obtain an element 
of $\boX(T)$ by the following construction. 
For all edges $[x,y]$ which lift to $[\widehat{x},\widehat{y}]$ we set 
$\phi_{x,y}=\phi(F(\widehat{x}),F(\widehat{y}))$
and for all triangles $[x,y,z]$ which lift to $[\widehat{x},\widehat{y},\widehat{z}]$ we 
define the $\Phi$ and $\delta$ invariants of $x,y,z$ as being equal to the corresponding 
invariants of the triple $\left(F(\widehat{x}),F(\widehat{y}),F(\widehat{z})\right)$. These 
data fit by construction as an element of $\boX(T)$. The two maps we have constructed are 
inverse one of the other. This ends the proof.

\end{proof}

\section{\label{resol}Solving the equations}

The aim of this part is to show how to construct solutions of the equations involved in $\boX(T)$ in a systematic way. 
The key lemma is the following:

\begin{lem}\label{resou3}
Let $m_{12}$, $m_{23}$, $m_{31}$ be three complex number different
from $0$ and $1$. From these numbers,  define
$\phi_{i,j}=|m_{ij}/(m_{ij}-1)|^2$ for all $i,j$. 
For any family $\left(\Phi_{ijk}\right)_{i,j,k}$ of complex numbers
satisfying the conditions
\begin{eqnarray}
\Phi_{ijk} & = & \Phi_{jki} \quad = \quad \overline{\Phi_{ikj}}\nonumber\\
\lvert\Phi_{ijk}\rvert & = & \sqrt{\phi_{ij}\phi_{jk}\phi_{kj}}\nonumber\\
\Delta_{ijk} & = & 1-\phi_{ij}-\phi_{jk}-\phi_{ki}+
\Phi_{ijk}+\overline{\Phi_{ijk}}<0\nonumber
\end{eqnarray}

 the following set of equations 

\begin{eqnarray}
 \delta^i_{jk}\delta^i_{kj} & = &\phi_{jk}\nonumber\\
 m^k_{ij}&=&\frac{1}{\Delta_{ijk}\phi_{ik}\phi_{jk}
 }(\phi_{ik}\phi_{jk}(\Phi_{ijk}-\phi_{ij})+\nonumber\\ 
&&\phi_{ik}(\phi_{ij}\phi_{jk}-\Phi_{ijk})\delta^i_{kj}+\phi_{jk}(\phi_{ij}\phi_{ik}-\Phi_{ijk})
\ba{\delta^j_{ki}}+\Phi_{ijk}(1-\phi_{ij})\delta^i_{kj}\ba{\delta^j_{ki}})\nonumber\\ 
0 & = &|\delta^i_{jk}|^2(1-\phi_{ik})+2\Re
  [\delta^i_{jk}(\Phi_{ikj}-\phi_{jk})]+\phi_{jk}(1-\phi_{ik}) \nonumber
\end{eqnarray}
have two distinct solutions in the variables $\delta^i_{jk}$ provided
that $\phi_{ij}\ne 1$ for all $i$ and $j$ in $\{1,2,3\}$. 
\end{lem}
\begin{proof}
Geometrically, the lemma has the following interpretation: let
$C_1,C_2,C_3$ be three complex lines in generic position. Their
position is parametrized by  the invariants $\phi_{ij}$ and
$\Phi_{ijk}$. The hypothesis on these invariants means that any two
complex lines are neither orthogonal nor asymptotic.  

The invariant $m_{ij}$ specifies a Lagrangian reflexion $I_{ij}$
swapping $C_i$ and $C_j$ in the following way: let $p_1$ and $p_2$ be
two points in $\partial C_1$ and $\partial C_2$ respectively such that
$m[(C_1,p_1),(C_2,p_2)]=m_{12}$. Then the second part of lemma \ref{elemRplan} 
tells us that there is a unique lagrangian involution swapping $C_1$ and $C_2$
and sending $p_1$ on $p_2$. This involution depends on $p_1$ and $p_2$
only through the data of $m_{12}$. In some sense, the involution
$I_{12}$ is the geometric realization of the invariant $m_{12}$.

A solution of the equations is equivalent to a triple of points
$p_1,p_2,p_3$ lying respectively in $\partial C_1,\partial C_2$ and
$\partial C_3$ such that for all $i,j$, $I_{ij}p_i=p_j$. Fixing a
reference complex line, say $C_1$, we see that a solution of the
equations is given by a fixed point of the product
$I_{31}I_{23}I_{12}$. This product is an anti-holomorphic isometry of
$C_1$ preserving the boundary, hence it has two distinct fixed points
on this circle. This proves the lemma.

 If some of the $\phi_{ij}$ are equal to one, then the corresponding
 complex lines are asymptotic. The same argument as above applies but
 the points $p_i$ may lie at the intersection of two complex lines
 which is not allowed in our settings. Hence, there are less than 2
 admissible solutions but there are still some degenerate ones.

\end{proof}

One can apply the preceding lemma for each triangle of a triangulation
at the same time. This is described in the following part.

Let $T$ be an ideal triangulation of
$\hat{\Sigma}_{g,p}$. We define a decoration space of
$T$ which is related to $\boX(T)$ but which is somewhat simpler: let
$\boM(T)$ be the set of triple $(\phi,\Phi,m)$ where: 
\begin{itemize} 
\item
$\phi$ and $m$ are functions defined on the set of oriented edges to $\C$ satisfying the following relations: 
$$\phi_{ij}=\left|\dfrac{m_{ij}}{m_{ij}-1}\right|^2 \mbox{ and } m_{ji}=\overline{m_{ij}}.$$
Note that the $\phi$ invariant is redundant as it is a function of $m$ but we keep it for the coherence of the notation.
\item
$\Phi$ is a $\C$ valued function defined on ordered faces of $T$
satisfying the following equations for all ordered faces $(i,j,k)$: 

$$\Phi_{ijk}=\Phi_{jki}=\overline{\Phi_{ikj}} \mbox{ and } \Delta_{ijk}<0$$
\end{itemize}

We denote by $\boX^{nd}(T)$ (resp. $\boM^{nd}(T)$) the non-degenerate part
of $\boX(T)$ (resp. $\boM(T)$) by which we mean the open set
of triples $(\phi,\Phi,\delta)$ such that $\phi_{ij}\ne1$ for all $i$
and $j$ (resp. the triples $(\phi,\Phi,m)$ such that
$\phi_{ij}\ne 1$ for all $i,j$). 

The following proposition is a direct consequence of the preceding lemma. 
\begin{prop}\label{cover}
The natural map $\boX^{nd}(T)\to\boM^{nd}(T)$ sending
$(\phi,\Phi,\delta)$ to $(\phi,\Phi,m)$ is a covering of order $2^{N}$
where $N$ is equal to the number of triangles in $T$. 
\end{prop}

This proposition explains than we can
solve the equations in a simple way: we fix arbitrarily the $\phi$ and
$m$ invariants, and then we solve (with a computer) the remaining
equations in $\delta$. The important point given by  the proposition
is that we are sure to obtain $2^N$ solutions in the non-degenerate
case. The simple structure of the map from $\boX^{nd}(T)$ to
$\boM^{nd}(T)$ should allow us to describe precisely the
representation space but it still does not seem to be an easy task and
we do not have done it yet. 
 
\section{Controlling the holonomy of the cusps\label{cusps}}
\subsection{The general case}
Consider a pair $(\rho,F)\in\mathfrak{R}_{g,p}$ and denote  as usual by $c_i$ the curve in $\Sigma_{g,p}$ enclosing $x_i$. Since
$\rho(c_i)$ stabilizes a flag $F_i=(C_i,p_i)$, it might  be either

\begin{itemize}
\item loxodromic, in which case its second fixed point belongs to $C_i$, 
\item parabolic, in which case $p_i$ is its unique fixed point,
\item a complex reflection, in which case its restriction to $C_i$ is
  the identity.
\end{itemize}

We wish to determine the type of  $\rho(c_i)$ in terms of the invariants $\phi$, $\Phi$ and $\delta$.
The loop $c_i$ encloses the vertex point $x_i$.
It may be written $c_i = \gamma \nu\gamma^{-1}$, where $\gamma$ is a path 
connecting the base point to one of the vertices of the hexagonation which is adjacent to the point $x_i$, and $\nu$ is       
a loop around $x_i$ which is composed of a succession of edges of the
the hexagonation connecting two edges of the triangulation. As a
consequence $\rho(c_i)$, may be written $M_\gamma NM_{\gamma}^{-1}$,
where $N$ is a product of elementary matrices which are all of
transfer type (see proposition \ref{isomflag}). Write 
$$N=\T_k \ldots \T_j\ldots \T_1,$$
where $\T_j$ is a matrix of transfer type:

$$\T_j=\begin{bmatrix} \mu_j & 0 & it_j\mu_j\\ 0 & \bar{\mu}_j/\mu_j & 0\\
0 & 0 & 1/\bar{\mu}_j\end{bmatrix},$$
\n and the $\mu_j$'s and $t_j$'s are written in terms of the invariants $\phi$,
$\Phi$ and $\delta$ as in proposition \ref{isomflag}. Computing the product, we obtain

$$N=\begin{bmatrix} \mu & 0 & K \\ 0 & \bar{\mu}/\mu & 0\\
0 & 0 & 1/\bar{\mu}\end{bmatrix},$$

\n where $\mu=\prod{\mu_j}$, and 

$$K = i\sum_{j=1}^{k}t_j\dfrac{\prod_{l=j}^{k}\mu_l}{\prod_{l=1}^{j-1}\bar\mu_l}.$$

\n We obtain thus that 

\begin{itemize}
\item $N$ is loxodromic if and only if $|\mu|\neq 1$,
\item $N$ is parabolic if and only if $|\mu|=1$ and $K\neq 0$, 
\item $N$ is a complex reflection if and only if $|\mu|=1$ and $K=0$. 
\end{itemize}

\subsection{Type preserving representations of the 1-punctured torus}

In this section, we focus on the special case of $\Sigma_{1,1}$, the 1-punctured torus. 
We first summarize the existing results about this case. We denote the fundamental group  of
$\Sigma_{1,1}$ by
$$\pi_{1,1} = \la a,b,c\,|\,\lbrack a,b\rbrack\cdot c=1\ra.$$
Recall that a representation of $\pi_{1,1}$ is said to be \textit{type preserving} if and only if $\rho(c)$ 
is a parabolic isometry. Note that there are two main types of parabolic isometries (see \cite{Go,Wi3} for more details):
\begin{enumerate}
\item\label{vert}  screw parabolic isometries. These parabolic elements preserve a complex, and thus a flag.
\item \label{hor}horizontal parabolic isometries, which preserves an $\R$-plane containing their fixed point. These isometries 
are also called non-vertical Heisenberg translations. No such parabolic isometry appears within the frame of the 
present work.
\end{enumerate}

\n As a consequence, we can separate the type-preserving representations of $\pi_{1,1}$ into two types, 
according to whether $\rho(c)$ is screw parabolic or horizontal parabolic.

\begin{enumerate}
\item If $\rho(c)$ preserves a complex line, then it is in the frame of this work. 
In this case, all the examples known  of a discrete, faithful and type preserving
representation are obtained by passing to an index 6 subgroup in a discrete, faithful
and type preserving representation of the modular group PSL(2,$\mathbb{Z}$).
 The latter representations have all been described by Falbel and Parker in \cite{FJP}. 
 This family of examples consists up to PU(2,1) of 6 topological components, 4 of which are
 points, and the two other are segments.
\item If $\rho(c)$ does not preserve a complex line, then it is a consequence of \cite{Wi4} that there exists a unique 
triple of Lagrangian reflections $(I_1,I_2,I_3)$ such that $\rho(a)=I_1I_2$ and $\rho(b)=I_3I_2$.
 In \cite{Wi3,Wi2}, all these type preserving representations of $\pi_{1,1}$ are described, and,
  among them, a 3-dimensional family of discrete and faithful representations is identified.
\end{enumerate}

\n We will now give  necessary and sufficient conditions for $\rho(c)$ to be parabolic, 
in the case it preserves a flag.\\
 
\begin{figure}
\begin{center}
\vspace{1cm}
\begin{pspicture}(-2,-2)(4,4)
\includegraphics[width=5cm,height=5cm,scale=0.4]{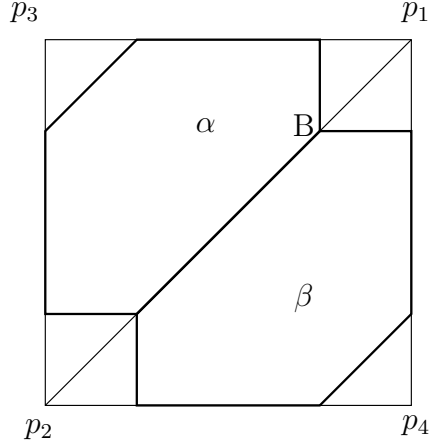}
\put(-5.2,-0.3){$p_2$}
\rput(0,-0.2){$p_4$}
\rput(0,5.3){$p_1$}
\rput(-5.2,5.3){$p_3$}
\rput(-1.5,3.8){B}
\rput(-2.8,3.8){$\alpha$}
\rput(-1.5,1.5){$\beta$}
\end{pspicture}
\caption[texte]{\protect\parbox[t]{12cm}{Ideal triangulation and hexagonation of a 1-punctured torus. 
The opposite sides of the square are identified.\label{torep}}}

\end{center}
\end{figure}

\n A triangulation of a 1-punctured torus is made of two triangles, which we call $\alpha$ and $\beta$ as
on figure \ref{torep}. We label the vertices as on figure \ref{torep}. The
decoration of this triangulation is given by :
\begin{itemize} 
\item Triangle  $\alpha$ $(p_1,p_2,p_3)$:  $\phi_{12}$, $\phi_{23}$, $\phi_{31}$, $\Phi_{123}$, $\delta^{1}_{23}$, 
$\delta^2_{31}$ and $\delta^3_{12}$.
\item Triangle $\beta$  $(p_1,p_2,p_4)$ :  $\phi_{12}$, $\phi_{24}$, $\phi_{41}$, $\Phi_{124}$, $\delta^{1}_{24}$,
$\delta^2_{41}$ and $\delta^4_{12}$.
\end{itemize}
 
 \n Note that because of the identification of the opposite sides of the square, the following relations hold:
 $$\phi_{23}=\phi_{14}\mbox{ and } \phi_{13}=\phi_{24}.$$

 We choose as a basepoint the vertex of the hexagonation which marked by B on figure \ref{torep}.
 Let us call $\nu^{i}_{jk}$ the oriented edge of the hexagonation turning around the vertex $p_i$ from the edge 
 $p_ip_j$ edge to the edge $p_ip_k$. As an example, $\nu^1_{24}$ is the oriented segment starting from 
 the point B (see figure \ref{torep}) and connecting the diagonal to the vertical side $p_1p_4$. The homotopy class  $c$
 is represented by the following sequence of edges, $ \nu^{1}_{23}\nu^{4}_{21}\nu^{2}_{31}\nu^{2}_{14}\nu^{3}_{12}
 \nu^{1}_{42}$, to which correspond the product of transfer type matrices $\T=\T^1_{42}\T^3_{12}\T^2_{14}\T^2_{31}
 \T^4_{21}\T^1_{23}$. Denote by $\mu^i_{jk}$ and $t^i_{jk}$ the two parameters in the matrix $\T^i_{jk}$ given by proposition
 \ref{isomflag}. The matrix $\T$ is upper triangular, and according to proposition \ref{isomflag}, its top left
 coefficient is 
 \begin{equation}\mu =
 \Theta\left( \dfrac{\delta^{1}_{23}\phi_{13}}{\Phi_{123}}\dfrac{\delta^{4}_{21}\phi_{14}}{\Phi_{421}}
 \dfrac{\delta^{2}_{31}\phi_{12}}{\Phi_{231}}\dfrac{\delta^{2}_{14}\phi_{24}}{\Phi_{214}}\dfrac{\delta^{3}_{12}\phi_{32}}
 {\Phi_{312}}\dfrac{\delta^{1}_{42}\phi_{12}}{\Phi_{142}}\right).
 \end{equation}

 We simplify this relation using the relations between the invariants ($\phi_{ij}=\phi_{ji}$, 
 $|\Phi_{ijk}|^2=\phi_{ij}\phi_{jk}\phi_{ki}$, and $\delta^{i}_{jk}\delta^{i}_{ki}=\phi_{ij}$). 
 This yields:
 
 $$|\mu|^2=\dfrac{|\delta^1_{23}\delta^2_{31}\delta^3_{12}|^2}{|\delta^1_{42}\delta^2_{14}\delta^4_{12}|^2}.$$
 
 We obtain as a direct consequence the following
 \begin{prop}
 Let $(\phi,\Phi,\delta)$ be a decorated triangulation of the punctured torus. The holonomy of a loop around
  the puncture is parabolic or a complex reflexion if and only if 
 \begin{equation}|\delta^1_{23}\delta^2_{31}\delta^3_{12}|=|\delta^1_{42}\delta^2_{14}\delta^4_{12}|
 \label{parabrefl}\end{equation}
 Moreover, the representation associated to the decoration is type preserving if and only if the relation
 (\ref{parabrefl}) is satisfied and the following relation holds
\begin{eqnarray}\label{nouveauK}
0 & \neq &\mu^1_{42}\mu^3_{12}\mu^2_{14}\mu^2_{31}\mu^4_{21}\mu^1_{23}\,t^1_{23}+
\dfrac{\mu^1_{42}\mu^3_{12}\mu^2_{14}\mu^2_{31}\mu^4_{21}}{\mu^1_{23}}\,t^4_{21}+
\dfrac{\mu^1_{42}\mu^3_{12}\mu^2_{14}\mu^2_{31}}{\mu^4_{21}\mu^1_{23}}\,t^2_{31}\nonumber\\
 & &+\dfrac{\mu^1_{42}\mu^3_{12}\mu^2_{14}}{\mu^2_{31}\mu^4_{21}\mu^1_{23}}\,t^2_{14}+
\dfrac{\mu^1_{42}\mu^3_{12}}{\mu^2_{14}\mu^2_{31}\mu^4_{21}\mu^1_{23}}\,t^3_{12}+
\dfrac{\mu^1_{42}}{\mu^3_{12}\mu^2_{14}\mu^2_{31}\mu^4_{21}\mu^1_{23}}\,t^1_{42}
\end{eqnarray}
  \end{prop}
 The relation (\ref{nouveauK}) is just an explicit version of $K\neq 0$, with $K$ as in
 the previous section.

\end{document}